\documentclass[11pt,twoside,reqno]{amsart}

\usepackage[
  a4paper,
  left=35mm,
  right=35mm,
  top=30mm,
  bottom=30mm,
  marginpar=25mm
]{geometry}

\usepackage[T1]{fontenc}
\usepackage{lmodern,microtype}
\usepackage{amsmath,amssymb,mathrsfs,stmaryrd,esint,accents}
\usepackage{booktabs}
\usepackage{enumerate}
\usepackage{enumitem}
\usepackage{comment,dirtytalk}
\usepackage{xcolor}
\usepackage{tikz,tikz-cd}
\usetikzlibrary{arrows.meta,calc}
\usepackage[unicode]{hyperref}
\usepackage{orcidlink}

\allowdisplaybreaks

\definecolor{darkcyan}{cmyk}{1,0,0,0.6}

\hypersetup{
  colorlinks=true,
  linkcolor=darkcyan,
  urlcolor=darkcyan,
  citecolor=darkcyan,
  pdftitle={Rough Vector Fields with Vanishing Lie Brackets},
  pdfauthor={Stefano Bianchini and Luca Talamini}
}

\makeatletter

\def\@settitle{%
  \begin{center}
    \normalfont\bfseries\LARGE
    \baselineskip=20pt\relax
    \@title\par
  \end{center}%
}

\def\@setauthors{%
  \begingroup
  \def\thanks{\protect\thanks@warning}%
  \trivlist
  \centering\small
  \@topsep30\p@\relax
  \advance\@topsep by -\baselineskip
  \item\relax
  \author@andify\authors
  \def\\{\protect\linebreak}%
  \textsc{\authors}%
  \endtrivlist
  \endgroup
}

\renewenvironment{abstract}{%
  \global\setbox\abstractbox=\vtop\bgroup
  \normalfont
  \vskip 1pc
  {\centering\normalsize\bfseries Abstract\par}
  \vskip 1pc
  \list{}{%
    \labelwidth\z@
    \leftmargin3pc
    \rightmargin\leftmargin
    \listparindent\normalparindent
    \itemindent\normalparindent
    \parsep\z@
    
  }%
  \item[]\normalfont\small
}{%
  \endlist\egroup
  \ifx\@setabstract\relax
    \@setabstracta
  \fi
}

\def\section{\@startsection{section}{1}%
  \z@
  {1.2\linespacing\@plus\linespacing}%
  {.5\linespacing}%
  {\normalfont\bfseries\centering}%
}

\def\subsection{\@startsection{subsection}{2}%
  {\parindent}%
  {.5\linespacing}%
  {-.5em}%
  {\normalfont\itshape}%
}

\def\subsubsection{\@startsection{subsubsection}{3}%
  {\parindent}%
  {.5\linespacing}%
  {-.5em}%
  {\normalfont\itshape}%
}

\def\ps@headings{%
  \ps@empty
  \def\@evenhead{%
    \normalfont\small\scshape
    \rlap{\thepage}%
    \hfil
    Stefano Bianchini and Luca Talamini%
    \hfil
  }%
  \def\@oddhead{%
    \normalfont\small\scshape
    \hfil
    Rough Vector Fields with Vanishing Lie Brackets%
    \hfil
    \llap{\thepage}%
  }%
}

\makeatother

\newtheorem{thm}{Theorem}[section]
\newtheorem{prop}[thm]{Proposition}

\newtheorem{lemma}[thm]{Lemma}

\theoremstyle{definition}

\theoremstyle{remark}

\numberwithin{equation}{section}

\newcommand*{\dif}{\mathop{}\!\mathrm{d}}

\newcommand{\mr}[1]{\mathrm{#1}}
\newcommand{\bs}[1]{\boldsymbol{#1}}

\newcommand{\wt}[1]{\widetilde{#1}}

\title[Rough vector fields with vanishing Lie brackets]%
{Rough vector fields with vanishing Lie  brackets}

\author[Stefano Bianchini]{%
  Stefano Bianchini\,\orcidlink{0000-0001-7886-3285}%
}
\address{Mathematics Area, SISSA,  Via Bonomea 265, 34136 Trieste, Italy}
\email{bianchin@sissa.it}

\author[Luca Talamini]{%
  Luca Talamini\,\orcidlink{0009-0008-5147-5184}%
}
\address{Mathematics Area, SISSA, Via Bonomea 265, 34136 Trieste, Italy}
\email{ltalamin@sissa.it}

\subjclass[2020]{35F35,37C10,46E35}

\keywords{Sobolev vector fields, Regular Lagrangian flows, Lie brackets}

\date{}

\begin{document}

\begin{abstract}
We consider two vector fields $X, Y \in W^{1,p}_{loc}(\mathbb R^d) \cap \mathbf L^\infty$ with bounded divergence. 
If $p \geq 2$, this paper shows that $DX \,Y - DY \, X  = 0$ almost everywhere if and only if their Regular Lagrangian Flows commute. If instead $1 \leq p < 2$, we present an example of two divergence free locally Sobolev vector fields in $\mathbb R^2$ with vanishing Lie brackets whose flows do not commute.  In the regime $2 \leq p  < \infty$ this paper also yields a quantitative commutation result not depending on the Lipschitz constant for  vector fields with $[X, Y] \neq 0$, which is new even in the smooth setting. 
\end{abstract}

\maketitle

\section{Introduction}
 If $X \in W^{1,p}_{loc}(\mathbb R^d)\cap \mathbf L^\infty$ is a Sobolev vector field with $\mathrm{div}\, X \in \mathbf L^\infty$, DiPerna and Lions \cite{DiPernaLions} (and then Ambrosio \cite{Ambrosio} in the $BV$ setting) proved that there is a unique natural flow map $\Phi_X^t: \mathbb R^d \to \mathbb R^d$ such that for all $T >0$ there is $C >0$ such that there holds
\begin{equation}\label{eq:comprb}
    (\Phi_X^t)_\sharp \mathcal L^d  \leq C \, \mathscr{L}^d \qquad \text{ for all $|t| \leq T$}
\end{equation}
and for almost every $x$, the map $t \mapsto \Phi_X^t(x)$ is an absolutely continuous solution to 
\begin{equation}\label{eq:traj}
\begin{cases}
& \frac{\dif}{\dif t} \Phi_X^t(x) = X(\Phi_X^t(x))\\
& \Phi_X^0(x) = x.
\end{cases}
\end{equation}
The map $\Phi_X^t$ is called the \emph{Regular Lagrangian flow} of the vector field $X$. 
The classical Frobenius theorem states that the (classical) flows $\Phi_X^t, \Phi_Y^s$ of two smooth vector fields $X, Y : \mathbb R^d \to \mathbb R^d$ commute
$$
\Phi^t_X \circ \Phi^s_Y = \Phi^s_Y \circ \Phi^t_X, \qquad \forall \; t, s  \geq 0
$$
if and only if their Lie bracket vanishes, i.e.
\begin{equation}\label{eq:LB}
    [X, Y] \,:= \,DX \, Y -DY \, X = 0.
\end{equation}
The aim of this paper is to extend the classical Frobenius theorem to the setting of Regular Lagrangian Flows.

Let $X, Y \in W^{1,p}$ be two bounded Sobolev fields, with bounded divergence. If their Regular Lagrangian Flows commute, it is not difficult (see e.g. \cite{ColomboTione}) to show that $[X, Y] = 0$ almost everywhere, therefore this paper and the previous literature on the topic only focus on the opposite implication. Assuming that \eqref{eq:LB} is satisfied almost everywhere, the following results are known about the commutation of their Regular Lagrangian flows.
\begin{enumerate}
    \item A first paper \cite{RampazzoSussmann} proves that the flows commute when $X, Y$ are locally Lipschitz.
    \item Later, \cite{ColomboTione} proved that the flows commute assuming additionally that for every $t\in\mathbb R$, the flow map $\Phi_X^t$ is weakly Lie differentiable in the direction $Y$, with locally bounded derivative.
    \item The paper \cite{RigoniStepanovTrevisan} proved the result when one of the vector fields is Lipschitz and the other belongs to $W^{1,p}$ and admits a unique Regular Lagrangian flow.
    \item Finally \cite{RebucciZizza} proved the commutativity of the flows for bounded, continuous Sobolev vector fields with bounded divergence on $\mathbb R^2$.
\end{enumerate}

The main result of this paper is as follows. Here $d$ is any dimension $d \geq 2$.
\begin{thm}\label{thm:MT}
   Let $X, Y \in W^{1,p}_{loc}(\mathbb R^d) \cap \mathbf L^\infty$ with $\mathrm{div}\, X$, $\mathrm{div}\, Y \in \mathbf L^\infty$, and assume that $[X,Y] = 0$. If $p \geq 2$, then for all $t, s > 0$, 
    \begin{equation}\label{eq:CF}
    \Phi^t_X \circ \Phi^s_Y(x)  = \Phi^s_Y\circ  \Phi^t_X(x) \qquad \text{for a.e. $x$}.
    \end{equation}
    If instead $1 \leq p < 2$, there are $X, Y \in W^{1,p}_{loc}(\mathbb R^d) \cap \mathbf L^\infty$ with zero divergence, such that \eqref{eq:CF} fails for $x$ in an open set with positive measure.
\end{thm}

Some comments are in order. 
\subsection{Lie differential and commutation} To understand why a different approach is needed compared to the existing literature,  assume that the flow $\Phi^t_X$ is Lie-differentiable in the direction of $Y$, with locally bounded derivative (see \cite{ColomboTione}). For the purpose of this sketch, assume that there exists a function $\vec \rho_t \in \mathbf L^\infty([0,T] \times \mathbb R^{d}; \, \mathbb R^d)$ such that
    $$
    \frac{1}{h} \Big( \Phi^t_X\circ \Phi^h_Y- \Phi^t_X \Big) \, \longrightarrow \, \vec \rho_t \qquad \text{in $\mathbf L^1_{loc}$ as $h \to 0^+$}.
    $$
    At least formally $\vec \rho_t$ satisfies 
    \begin{equation}\label{eq:Lieder}
    \partial_t \vec \rho_t = DX(\Phi^t_X)\, \vec \rho_t, \qquad  \vec \rho_0 = Y.
     \end{equation}
    Notice that by applying Gronwall's estimate along almost every trajectory of $X$, the above equation has a unique solution $\vec \rho_t$ in $\mathbf L^\infty$. However, due to the assumption $[X, Y] = 0$, it is immediate to check that $Y({\Phi}^t_X(x))$ is, in fact, another bounded solution of \eqref{eq:Lieder}, and therefore 
    $$
   \vec  \rho_t = Y(\Phi^t_X).
    $$
    Therefore 
    $$
    \frac{\dif}{\dif s} \Phi^t_X\circ\Phi^s_Y(x) = Y(\Phi^t_X\circ\Phi^s_Y(x))
    $$
    so that $s \mapsto \Phi^t_X\circ\Phi^s_Y(x)$ is a trajectory of the flow of $Y$ for a.e. $x$, starting at the point $\Phi^t_X(x)$. Since it satisfies \eqref{eq:comprb} and $\eqref{eq:traj}$, again by uniqueness of the DiPerna Lions flows one must have that it is equal to $\Phi^s_Y \circ \Phi^t_X(x)$ at a.e. $x$.
    
    Results about the differentiability of the Regular Lagrangian Flow of Sobolev fields are available in two dimensions under continuity assumptions \cite{Marconi} (and, classically, for Lipschitz vector fields in general dimension), but it is well known, see \cite{Marconi,Jabin} that the RLF of Sobolev vector fields need not  be differentiable or inherit any Sobolev regularity. In particular, there seems not to be an easy way to fix the previous argument in a general Sobolev setting.

\subsection{Quantitative commutation}
For $p > 2$, strictly greater than the critical exponent, the proof also provides a quantitative commutation result for vector fields in $W^{1,p}$, as follows.
\begin{thm}\label{thm:QC}
    Let $p > 2$ and let $X, Y \in W^{1,p}_{loc}(\mathbb R^d) \cap \mathbf L^\infty$ with $\mathrm{div}\, X$, $\mathrm{div}\, Y \in \mathbf L^\infty$. Then  
    \begin{equation}\label{eq:QC}
        \sup_{|s|, |t| \leq T} \int_{B_r} \left|\Phi^t_X \circ \Phi^s_Y(x) - \Phi^s_Y \circ \Phi^t_X(x)\right| \dif x \lesssim  \left[ \log \left( 1 + \frac{1}{\|[X, Y] \|_{\mathbf L^1(B_R)}}\right) \right]^{-p}
    \end{equation}
    where $$R = r +  T(\|X\|_\infty + \|Y\|_\infty).$$ The implicit constant depends on $d, p, T, r$ and the norms of $X, Y$ in $B_{R+1}$.
\end{thm}

In  the critical case $X, Y \in W^{1,2}$, one readily deduces the existence of an analogous quantitative result, using the compactness and stability properties of the Regular Lagrangian Flows. However, the uniform modulus that we can obtain when $p = 2$ is not as clean as in the case $p > 2$, as it will be clear from the proof below.

\vspace{0.3cm}
    Since $X, Y \in \mathbf L^\infty$, by a standard localization argument it is easy to see that it is equivalent to prove the commutativity results for vector fields in $W^{1,p}(\mathbb R^d)$. To avoid technical issues, we will restrict to this class from now on.

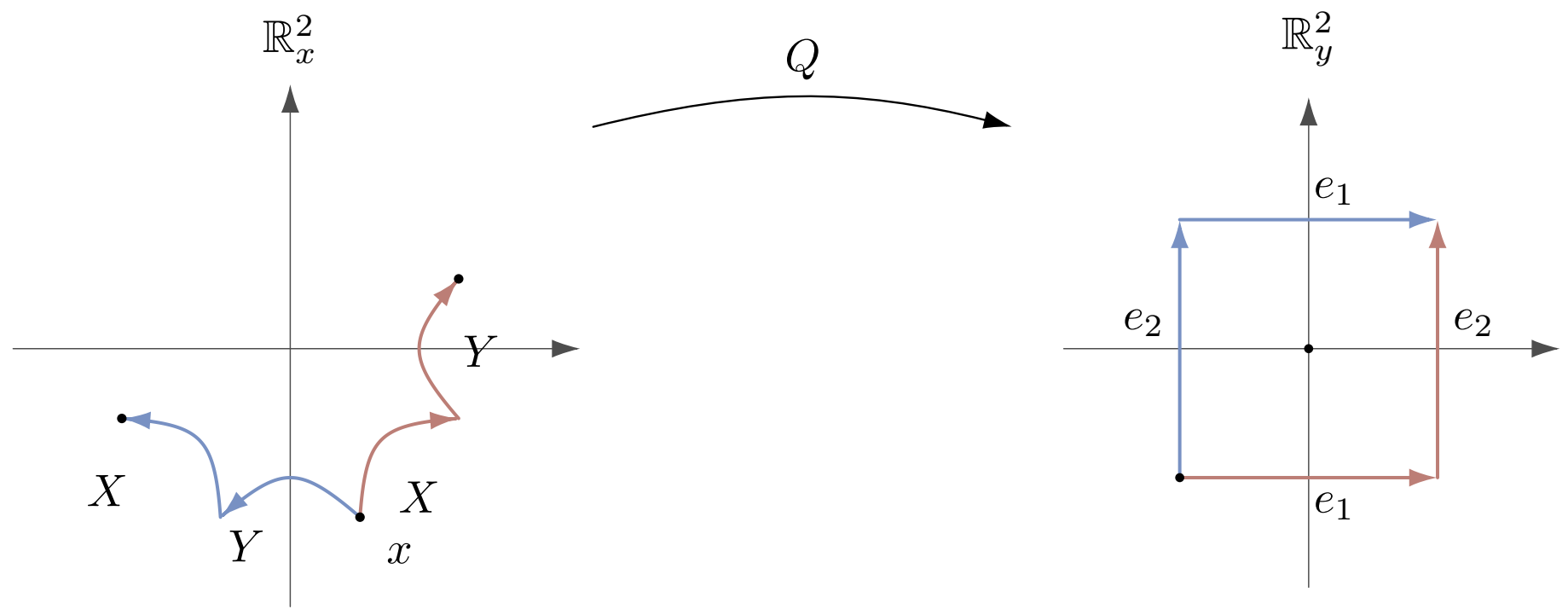
\begin{figure}
    \centering

\definecolor{trajectoryRed}{RGB}{188,126,118}
\definecolor{trajectoryBlue}{RGB}{120,145,195}

\begin{tikzpicture}[
  line cap=round,
  line join=round,
  >={Latex[length=2.2mm,width=1.4mm]},
  axis/.style={draw=black!70, line width=0.28pt, ->},
  maparrow/.style={draw=black, line width=0.45pt, ->},
  redpath/.style={draw=trajectoryRed, line width=0.75pt, ->},
  bluepath/.style={draw=trajectoryBlue, line width=0.75pt, ->},
  lab/.style={font=\small}
]

\begin{scope}

\draw[axis] (-2.15,0) -- (2.25,0);
\draw[axis] (0,-2.0) -- (0,2.05);
\node[lab] at (0,2.40) {$\mathbb{R}^2_x$};

\draw[redpath]
  plot[domain=0:1,samples=80,variable=\u]
  ({sqrt((2*\u-1)^2+1)*cos(atan2(-1,2*\u-1)/2)},
   {sqrt((2*\u-1)^2+1)*sin(atan2(-1,2*\u-1)/2)});

\draw[redpath]
  plot[domain=0:1,samples=80,variable=\u]
  ({sqrt(1+(2*\u-1)^2)*cos(atan2(2*\u-1,1)/2)},
   {sqrt(1+(2*\u-1)^2)*sin(atan2(2*\u-1,1)/2)});

\draw[bluepath]
  plot[domain=0:1,samples=80,variable=\u]
  ({sqrt(1+(2*\u-1)^2)*cos((-180+atan2(1-2*\u,1))/2)},
   {sqrt(1+(2*\u-1)^2)*sin((-180+atan2(1-2*\u,1))/2)});

\draw[bluepath]
  plot[domain=0:1,samples=80,variable=\u]
  ({sqrt((2*\u-1)^2+1)*cos((-270-atan2(2*\u-1,1))/2)},
   {sqrt((2*\u-1)^2+1)*sin((-270-atan2(2*\u-1,1))/2)});

\coordinate (x0) at
  ({sqrt(2)*cos(-67.5)},{sqrt(2)*sin(-67.5)});
\coordinate (xr) at
  ({sqrt(2)*cos(22.5)},{sqrt(2)*sin(22.5)});
\coordinate (xb) at
  ({sqrt(2)*cos(-157.5)},{sqrt(2)*sin(-157.5)});

\fill (x0) circle (1.05pt);
\fill (xr) circle (1.05pt);
\fill (xb) circle (1.05pt);

\node[lab,below right=2pt] at (x0) {$x$};

\node[lab] at (1.00,-1.15) {$X$};
\node[lab] at (1.48,-0.02) {$Y$};
\node[lab] at (-0.34,-1.53) {$Y$};
\node[lab] at (-1.42,-1.10) {$X$};

\end{scope}

\draw[maparrow] (2.35,1.72)
  to[bend left=14] node[above,lab] {$Q$} (5.60,1.72);

\begin{scope}[xshift=7.9cm]

\draw[axis] (-1.9,0) -- (1.95,0);
\draw[axis] (0,-1.85) -- (0,1.95);

\node[lab] at (0,2.40) {$\mathbb{R}^2_y$};

\draw[redpath] (-1,-1) -- (1,-1);
\draw[redpath] (1,-1) -- (1,1);

\draw[bluepath] (-1,-1) -- (-1,1);
\draw[bluepath] (-1,1) -- (1,1);

\fill (0,0) circle (1.0pt);
\fill (-1,-1) circle (1.0pt);

\node[lab] at (0.2,-1.22) {$e_1$};
\node[lab] at (1.28,0.2) {$e_2$};
\node[lab] at (-1.28,0.2) {$e_2$};
\node[lab] at (+0.2,1.22) {$e_1$};

\end{scope}
\end{tikzpicture}

    \caption{Example of two divergence free fields $X, Y \in W^{1,p}_{loc}(\mathbb R^2)$, $p < 2$, whose flows do not commute. They are constructed as pullbacks of $e_i$ through a branching map $Q: \mathbb R^2_x \to \mathbb R^2_y$. For example, in polar coordinates, one can take $Q(\rho, \theta) = (\rho, 2\theta)$.}
    \label{fig:nocomm}
\end{figure}

\subsection{The regime $p < 2$}  The example of non-commuting flows in the regime $p < 2$ is surprisingly simple, and it is constructed using a double cover of the plane $Q : \mathbb R^2_x \to \mathbb R^2_y$ (see Figure \ref{fig:nocomm}). The example was essentially already known from \cite{Nelson1959} (see also \cite{Fuglede1982, HeinonenKeith2011}), here we just need to check its Sobolev regularity, which is immediate. Say that, in polar coordinates
$$
Q(\rho, \theta) = (\rho, 2\theta).
$$ 
Vector fields with vanishing Lie brackets $X, Y$ producing non-commuting flows can be defined as the pullback of $e_1 = (1,0), e_2 = (0,1)$ via the map $Q$:
$$
X = (DQ)^{-1} \, e_1, \qquad Y = (DQ)^{-1} \, e_2.
$$
It is easy to see that they satisfy all the conditions of Theorem \ref{thm:MT} if and only if $p < 2$. Their flows do not commute: the reason is that, while clearly $e_1, e_2$ commute in $\mathbb R^2_y$, when they are pulled back via the branching map $Q$, the points $\Phi^t_X \circ \Phi^s_Y$ and $\Phi^s_Y \circ \Phi^t_X$ may well end up in different branches of the map $Q$ in the original space $\mathbb R^2_x$. The exponent $p = 2$ arises as a natural threshold also in \cite{HeinonenKeith2011} for the same geometric reasons.

\section*{AI use}

We started working on the problem in early 2025, but eventually got stuck and moved away from the topic. In September 2026, we were very surprised when GPT-6 Astra suggested (in a single response) the counterexample in the regime $p<2$. Moreover, this time after extensive interaction with GPT-6 Astra, in late September 2026 we obtained a first proof, completely different from the one presented here, that worked for $p>2$ in arbitrary dimension. Finally, in late October 2026, through further interactions with the same model, we obtained the present proof, which we then decided to rework and write up in its current form. The paper was written entirely by the authors.

\section{Short time estimate}

We start with some preliminary results, needed for the following.

\subsection{Maximal functions and truncations}\label{sec:preli} Recall that if $f \in W^{1,p}(\mathbb R^d; \, \mathbb R^m)$  then the \emph{maximal function} of $Df$ is defined as 
$$
g_f(x) := \mathcal{M}(|D f|)(x) \; =\;  \sup_{\rho > 0} \frac{1}{|B_\rho(x)|} \int_{B_\rho(x)} |Df(y)| \dif y.
$$
For the following statements, we refer to \cite{Stein1970}.
If $ p > 1$,  one has $$\|g_f\|_{\mathbf L^p} \; \leq\;  C_{d,p} \,\|f\|_{W^{1,p}}.$$
Then there exists a null set $N$ such that for all $x, y \in \mathbb R^d\setminus N$, there holds
\[
|f(x)-f(y)|
\; \le\;  C_d \; (g_f(x)+g_f(y))|x-y|.
\]
On the set 
$$G_\lambda:= \left\{x \in \mathbb R^d\; | \; g_f(x) \leq \lambda/(2 C_d)\right\}\setminus N,$$
the function $f$ is $\lambda$-Lipschitz, therefore, by the Kirszbraun Lipschitz extension theorem $f_{\mid G_\lambda}$ admits a Lipschitz extension $f_\lambda$,  with
$$
\|f_\lambda\|_{\infty} \leq \|f\|_{\infty}, \qquad \mr{Lip} \, f_{\lambda} \leq \lambda, \qquad f_\lambda(x) = f(x) \quad \text{for a.e.} \; x \in G_\lambda.
$$
\subsection{Lusin-Lipschitz estimates} It is well known that the Regular Lagrangian Flow of Sobolev vector fields satisfies Lipschitz estimates outside a small set, see \cite{CrippaDeLellis}. We recall here a version in our setting.
\begin{lemma}\label{lemma:Lusinlip}
    Let $p > 1$ and $X \in W^{1,p}(\mathbb R^d)$ be a vector field with bounded divergence, and let $g_X$ be its maximal function.   Then for every $\lambda > 0$ and $0 \leq h \leq T$ there is a set $\mathcal B_{\lambda,h}$ such that
    $$
    |\mathcal B_{\lambda,h}| \leq 2C \frac{ h}{\lambda}  \int_{\{g_X \geq \lambda/(2h)\}} g_X(x) \dif x
    $$
    and
    $$x \mapsto \Phi^h_X(x) \qquad \text{is $e^{2C_d \lambda}$-Lipschitz in $\mathbb R^d \setminus \mathcal B_{\lambda, h}$}.$$
\end{lemma}
\begin{proof}
    Up to a mollification argument, we can assume that $X$ is a smooth vector field. Then we obtain 
    $$
    \begin{aligned}
    \frac{\dif}{\dif t} |\Phi^t_X(x) - \Phi^t_X(y)| & \leq |X(\Phi^t_X(x))- X(\Phi^t_X(y))| \\
    & \leq\;  C_d \,(g_X(\Phi^t_X(x)) + g_X(\Phi^t_X(y)))|\Phi^t_X(x)- \Phi^t_X(y)|.
    \end{aligned}
    $$
    By Gronwall's Lemma, 
    $$
    |\Phi^h_X(x) - \Phi^h_X(y)| \leq |x-y| \exp\left( C_d\int_0^h(g_X(\Phi^t_X(x)) + g_X(\Phi^t_X(y))) \dif t\right).
    $$
    Define for $x \in \mathbb R^d$
$$
G_{\lambda, h}(x) := \int_0^h g_X(\Phi^t_X(x)) \dif t
$$
    and define the set 
    $$
    \mathcal B_{\lambda, h} := \left\{x \in \mathbb R^d \; | \; G_{\lambda, h}(x) \geq \lambda\right\}.
    $$
    Then for a.e. $x,y \in \mathbb R^d \setminus \mathcal B_{\lambda, h}$, there holds
    $$
    |\Phi^h_X(x) - \Phi^h_X(y)| \; \leq \; e^{(2C_d \lambda)} \, |x-y|.
    $$
    Moreover, 
$$
\begin{aligned}
    |\mathcal B_{\lambda, h}| & \leq \frac{1}{\lambda} \int_{\mathcal B_{\lambda, h}} G_{\lambda, h}(x) \dif x \\
    &  = \frac{1}{\lambda} \int_{\mathcal B_{\lambda, h}} \int_0^h g_X(\Phi^t_X(x)) \mathbf 1_{\{g_X(\Phi^t_X(x)) \leq \lambda/(2h)\}}(x) \dif t \dif x \\
    & \qquad + \frac{1}{\lambda} \int_{\mathcal B_{\lambda, h}} \int_0^h g_X(\Phi^t_X(x)) \mathbf 1_{\{g_X(\Phi^t_X(x)) > \lambda/(2h)\}}(x) \dif t \dif x \\
     & \leq \frac{1}{2}|\mathcal B_{\lambda, h}| +  \frac{Ch}{\lambda} \int_{\mathbb R^d} g_X(x) \mathbf 1_{\{g_X(x) \geq \lambda/(2h)\}}(x) \dif x.
\end{aligned}
$$
    where we have used \eqref{eq:comprb} to perform the change of variable in the last inequality. Since the first term in the last line can be absorbed in the left hand side, this concludes the proof. 
\end{proof}

\subsection{A short-time estimate}
	
The main goal of this section is to prove the following preliminary proposition, that is a short time estimate on the commutativity of flows when $t, s \sim h$.

\begin{prop}\label{prop:2}
   Let $X, Y \in W^{1,p}(\mathbb R^d) \cap \mathbf L^\infty$ be vector fields with bounded divergence. Fix $h > 0$ and $r > 0$. Then for all $0 \leq t, s \leq h$ there holds
   \begin{equation}
   \begin{aligned}
       \int_{B_r} \Big|\Phi_X^t \circ \Phi_Y^s(x) - \Phi_Y^s \circ \Phi_X^t(x) \Big| \dif x & \lesssim h \cdot \left|\left\{ g \geq h^{-1}\right\} \cap B_R \right|+ h^2 \,\big\| [X,Y]\big\|_{\mathbf L^1(B_R)}
          \end{aligned}
   \end{equation}
   where 
   \begin{equation}\label{eq:Rrh}
       R \equiv R_{r, h}:=  r + h(|X|_\infty + |Y|_\infty)
       \end{equation}
       and $g = g_X + g_Y$.
\end{prop}

The following is a preliminary estimate for smooth vector fields.
\begin{lemma}\label{lemma:lip}
    Let $X, Y: \mathbb R^d \to \mathbb R^d$ be bounded, Lipschitz vector fields with Lipschitz constant at most $\lambda$. Then there holds, for all $0 \leq t, s \leq h :=2C_d \lambda^{-1}$,
    \begin{equation}
        \int_{B_r}\Big|\Phi^t_{X} \circ \Phi_{Y}^s(x) - \Phi_Y^s \circ \Phi^t_{X}(x)\Big| \dif x  \; \lesssim_d  \,t \,s \,    e^{\lambda(s+t)}\, \big\| [X,Y]\big\|_{\mathbf L^1(B_{R_{r, h}})}.
    \end{equation}
    where the implicit constant is a dimensional constant.
\end{lemma}
\begin{proof}
Without loss of generality we prove the result for $t, s \geq 0$.  By a standard mollification argument, it is sufficient to prove the result for smooth vector fields. First consider the quantity 
$$
\mr E^t(z) := D \Phi_X^t(z)  Y(z) - Y(\Phi^t_X(z)), \qquad t \geq 0, \quad z \in \mathbb R^d
$$
Differentiating with respect to $t$, we obtain 
$$
\partial_t \mr E^t(z) = D X(\Phi^t_X(z))\mr E^t(z) + [X,Y](\Phi^t_X(z)), \qquad \mr E_0(z) = 0.
$$
Therefore, by Gronwall's Lemma we obtain 
\begin{equation}\label{eq:Etest}
|\mr E^t(z)| \leq e^{\lambda t} \int_0^t \big| \, [X,Y](\Phi^p_X(z)) \big| \dif p.
\end{equation}
Next, consider the function
$$
[0,s] \ni u \mapsto \mr C^{t,u}(x) := \Phi^t_{X} \circ \Phi_{Y}^u(x) - \Phi_Y^u \circ \Phi^t_{X}(x).
$$
Its derivative with respect to $u$ is
$$
\begin{aligned}
\partial_u \mr C^{t,u}(x) & =  D\Phi^t_X\big(\Phi_{Y}^u(x)\big) \,Y(\Phi^u_Y(x)) - Y(\Phi^u_Y \circ \Phi^t_X(x)) \\
& \qquad \qquad = \Big(D\Phi^t_X\big(\Phi_{Y}^u(x)\big) \,  Y(\Phi^u_Y(x)) - Y ( \Phi^t_X \circ \Phi^u_Y(x)) \Big) \\
& \qquad \qquad \qquad \qquad  + \Big(Y (\Phi^t_X \circ \Phi^u_Y(x)) - Y(\Phi^u_Y \circ \Phi^t_X(x)) \Big). 
\end{aligned}
$$
Therefore, by triangular inequality, using \eqref{eq:Etest}
$$
\begin{aligned}
    |\partial_u \mr C^{t,u}(x)| & \leq |\mr E^t(\Phi_Y^u(x))| + \lambda |\mr C^{t,u}(x)| \\
    & \leq e^{\lambda t} \int_0^t\big|\, [X,Y]\big(\Phi^p_X \circ \Phi_Y^u(x)\big)\, \big| \dif p + \lambda |\mr C^{t,u}(x)|
\end{aligned}
$$

Therefore again by Gronwall's Lemma we obtain 
$$
\begin{aligned}
    |\Phi^t_{X} \circ \Phi_{Y}^s - \Phi_Y^s \circ \Phi^t_{X}(x)| & \leq e^{\lambda(s+t)} \int_0^s  \int_0^t \big|\, [X,Y]\big(\Phi^p_X \circ \Phi_Y^u(x)\big)\, \big| \dif p \dif u
\end{aligned}
$$
A final integration in $x$ yields 
$$
\begin{aligned}
\int_{B_r} \big|\, \Phi^t_{X} \circ \Phi_{Y}^s - \Phi_Y^s \circ \Phi^t_{X}(x)\, \big| \dif x & \leq e^{\lambda(s+t)} \int_0^s \int_0^t \Bigg( \int_{B_r}  \big| \,[X,Y] (\Phi^p_X\circ \Phi_Y^u(x)\big)\big| \dif x \Bigg)\,  \dif p \dif u\\
& \leq C \,t s\,e^{\lambda(s+t)}\int_{B_R}\big| \, [X,Y](x) \, \big| \dif x
\end{aligned}
$$
where in the last inequality we use the compression bounds on the flows of $X, Y$ for small times $0 \leq t, s \leq 2C_d\lambda^{-1}$, which imply 
$$
(\Phi^p_X)_\sharp \mathcal L^d, \quad (\Phi^u_Y)_\sharp \mathcal L^d  \leq C \mathcal L^d
$$
 where $C$ is a dimensional constant.
\end{proof}

The following standard Lemma is an estimate of how much the flow of $X$ can differ from the flow of the truncated field $X_\lambda$.
\begin{lemma}\label{lemma:flowdist}
    If $X \in W^{1,p}(\mathbb R^d) \cap \mathbf L^\infty$ has bounded divergence, and $X_\lambda$ is its $\lambda$-Lipschitz truncation, then there holds for all $0 \leq t \leq T$
    $$
     \int_{B_r} \Big|\Phi^t_{X}(x) - \Phi^t_{X_\lambda}(x)\Big| \dif x \leq C \,t\, e^{\lambda t} \left| \left\{ X \neq X_\lambda\right\} \cap B_{r + t |X|_\infty}\right|.
    $$
\end{lemma}
\begin{proof}
    We have
$$
\begin{aligned}
\partial_t \Big( \Phi^t_{X}(x) -  \Phi^t_{X_\lambda}(x) \Big)  
 & = \Big(X(\Phi^t_X(x)) - X_\lambda(\Phi^t_{X}(x))\Big)  
 \\
 & + \,\Big( X_\lambda(\Phi^t_{X}(x)) - X_\lambda(\Phi^t_{X_\lambda}(x)) \Big)
\end{aligned}
$$
By estimating the second term using the Lipschitz regularity of $X_\lambda$, again by Gronwall's Lemma we obtain 
    $$
    |\Phi^t_{X}(x) - \Phi^t_{X_\lambda}(x)| \leq e^{\lambda t}\int_0^t |X(\Phi^s_X(x)) - X_\lambda(\Phi^s_{X}(x))| \dif s
    $$
    Integrating in $x$ 
    $$
    \begin{aligned}
    \int_{B_r} |\Phi^t_{X}(x) - \Phi^t_{X_\lambda}(x)| \dif x & \leq C\, t \, e^{\lambda t} \int_{B_R} |X(x) -X_\lambda(x)| \dif x\\
    & \leq \|X\|_\infty\, C\, t\, e^{\lambda t} |\{ X \neq X_\lambda\}\cap B_R|
    \end{aligned}
    $$
    where in the first inequality we used the compression bounds \eqref{eq:comprb} for the Regular Lagrangian flow of $X$.
\end{proof}

\begin{proof}[Proof of Proposition \ref{prop:2}]
   We set $$\lambda := 2C_d/h$$ and we  consider the truncated fields $X_\lambda, Y_\lambda$. We can write
    \begin{equation}\label{eq:tria}
    \begin{aligned}
      \int_{B_r} \Big|\Phi_X^t \circ \Phi_Y^s(x) - \Phi_Y^s   \circ \Phi_X^t(x)&  \Big| \dif x  \leq \int_{B_r} \Big|\Phi_{X_\lambda
      }^t \circ \Phi_{Y_\lambda}^s(x) - \Phi_{Y_\lambda}^s \circ \Phi_{X_\lambda}^t(x) \Big| \dif x  \\
      & \quad + \int_{B_r} \Big|\Phi_{X_\lambda
      }^t \circ \Phi_{Y_\lambda}^s(x) - \Phi_{X}^t \circ \Phi_{Y}^s(x) \Big| \dif x  \\
      & \quad + \int_{B_r} \Big|\Phi_{Y_\lambda
      }^s \circ \Phi_{X_\lambda}^t(x) - \Phi_{Y}^s \circ \Phi_{X}^t(x) \Big| \dif x.
    \end{aligned}
    \end{equation}
    Notice that by applying Lemma \ref{lemma:flowdist} twice, there holds
    $$
    \begin{aligned}
        \int_{B_r} \Big|\Phi_{X_\lambda
      }^t \circ \Phi_{Y_\lambda}^s(x) - \Phi_{X}^t & \circ \Phi_{Y}^s(x) \Big| \dif x   \lesssim  \int_{B_r} \Big|\Phi_{X_\lambda
      }^t \circ \Phi_{Y_\lambda}^s(x) - \Phi_{X_\lambda
      }^t \circ \Phi_{Y}^s(x)\Big|\dif x \\
      & \quad + \int_{B_r} \Big|\Phi_{X_\lambda
      }^t \circ \Phi_{Y}^s(x) - \Phi_{X
      }^t \circ \Phi_{Y}^s(x)\Big|\dif x \\
      & \lesssim\, e^{\lambda t} \int_{B_r}   \Big|\Phi^s_{Y_\lambda}(x) -\Phi^s_Y(x)\Big| \dif x \, +\,  t e^{\lambda t} \left|\{X \neq X_\lambda\}\cap B_R \right|\\
      & \lesssim \, s\, e^{\lambda t} e^{\lambda s} |\{ Y \neq Y_\lambda\}\cap B_R| \, +\,   t\, e^{\lambda t} |\{X \neq X_\lambda  \}\cap B_R|
    \end{aligned}
    $$
    The third term in \eqref{eq:tria} can be estimated in an entirely symmetrical way, while for the first one can use Lemma \ref{lemma:lip}. Combining all these estimates, and recalling the properties of the Lipschitz truncation of Section \ref{sec:preli},  one gets
     $$
    \begin{aligned}
        \int_{B_r} \Big|\Phi_{X
      }^t \circ \Phi_{Y}^s(x) - & \Phi_{Y}^s \circ \Phi_{X}^t(x) \Big| \dif x   \leq h^2 \int_{B_R} \big|\,[X_\lambda,Y_\lambda](x) \, \big| \dif x +h |B_R \cap\{g \geq h^{-1}\}|\\
      & \lesssim h^2 \int_{B_R} \big| \,[X,Y] \, \big| \dif x +  h |B_R \cap \{g \geq h^{-1}\}|
    \end{aligned}
    $$
    where we used that $|DX_\lambda|, |DY_\lambda| \leq \lambda = 2C_d/h$ to pass from the first to the second line and moreover that, on the good set $\{g \leq h^{-1}\}$, by locality of the weak gradient one has $DX = DX_\lambda$ and $DY = DY_\lambda$ almost everywhere.
    
\end{proof}

\section{A quantitative commutation estimate}
In this section we prove the first part of Theorem \ref{thm:MT}, i.e. that the flows commute whenever $p \geq 2$.

The following is the key proposition: the initial commutator is split into $4^k$ smaller commutators, up to removing a set with controlled measure. Here we recall that $g$ is the sum of the maximal functions of the derivatives of $X, Y$:
$$
g = g_X + g_Y
$$
as defined in Section \ref{sec:preli}.

\begin{prop}\label{prop:key}
Let $X,Y\in W^{1,p}(\mathbb R^d)\cap\mathbf L^\infty$,
with $p>1$ and bounded divergence. Let $0 \leq  s,t \leq T$,
and let $\{\lambda_j\}_{j\in\mathbb N}$ be a sequence of
positive numbers. Then, for every $r >0$ and $k>0$,
\begin{equation}\label{eq:key-estimate1}
\begin{aligned}
&\int_{B_r}
\Big|\Phi_X^t\circ\Phi_Y^s(x)
-\Phi_Y^s\circ\Phi_X^t(x)\Big|\dif x\\
&\quad\lesssim
\sum_{j=1}^{k}
\frac{2^j}{ \lambda_j}
\int_{\{g>\,2^j\lambda_j/(2T)\}}g(x)\dif x\\
&\qquad+
4^k\left(\prod_{j=1}^{k}e^{c_d\lambda_j}\right)
\int_{B_{R_k}}
\Big|
\Phi_X^{\frac{t}{2^k}}
\circ\Phi_Y^{\frac{s}{2^k}}(x)
-
\Phi_Y^{\frac{s}{2^k}}
\circ\Phi_X^{\frac{t}{2^k}}(x)
\Big|\dif x.
\end{aligned}
\end{equation}
where here 
$$
R_k = r + (1-2^{-k}) T(|X|_\infty +|Y|_\infty).
$$
and $c_d$ is a dimensional constant, and the implicit constant depends on $r, T, d$.
\end{prop}

\begin{proof}[Proof of Proposition \ref{prop:key}]

  {\bf 1.}  For convenience of notation, let us introduce, for all $j \in \mathbb N$,
    $$
    A_j := \Phi^{t/2^j}_X, \qquad B_j := \Phi^{s/2^j}_Y,
    $$
    $$
    \bs \Gamma_j := A_j \circ B_j - B_j\circ A_j.
    $$
   Notice the recursive relation
   \begin{equation}\label{eq:jdec}
   \begin{aligned}
    \bs \Gamma_j  = &  A_{j+1}  \circ \bs \Gamma_{j+1} \circ B_{j+1} \\
   & \quad +A_{j+1}  \circ B_{j+1}   \circ\bs \Gamma_{j+1} \\
   &  \qquad +  \bs \Gamma_{j+1} \circ B_{j+1} \circ A_{j+1} \\
   & \qquad\quad+ B_{j+1}   \circ\bs \Gamma_{j+1} \circ A_{j+1}.
   \end{aligned}
    \end{equation}
    Here we use the convention that for two maps $D_1, D_2 : \mathbb R^d \to \mathbb R^d$,
$$
D_1 \circ \bs \Gamma_j \circ D_2 \equiv D_1 \circ A_j \circ B_j\circ  D_2 - D_1 \circ B_j \circ A_j \circ D_2. 
$$
    
   Introduce the set of indices $\mathcal I^k := \{1,2,3,4\}^k$. An element $\alpha \in \{1,2,3,4\}^k$ is a vector of the form $$\alpha = (\alpha^1, \ldots, \alpha^k), \qquad \alpha^i \in \{1, 2, 3, 4\}.$$ Expanding the first commutator $\bs \Gamma_0$ recursively for $k \geq 1$ iterations, we therefore obtain for every $k \geq 1$ that 
   \begin{equation}
       \bs \Gamma_0 = \sum_{\alpha \in \mathcal I^k} (\mathcal L^1_{\alpha^1} \circ \ldots \circ \mathcal L^k_{\alpha^k}) \circ \bs \Gamma_k\circ  (\mathcal R^k_{\alpha^k} \circ \ldots \circ \mathcal R^1_{\alpha^1}).
   \end{equation}
   In particular, at step $k$, we decomposed $\bs \Gamma_0$ into the sum of $4^k$ commutators $\bs \Gamma_k$, left and right composed with some operators. Here $\mathcal L^j_{\alpha^j}$ is one of the four elements appearing in the left compositions in \eqref{eq:jdec}, more precisely:
   \begin{equation}\label{eq:left}
   \mathcal L^j_{\alpha^j} \in \Big\{A_{j}, \, A_j \circ B_j, \, \mr{Id}, \, B_{j} \Big\}.
      \end{equation}
   In an entirely similar way
   \begin{equation}\label{eq:right}
   \mathcal R^j_{\alpha^j} \in \Big\{B_j, \, \mr{Id}, \, B_j \circ A_j, \, A_j \Big\}.
\end{equation}
   Let us name 
   $$
   \wt{\mathcal L}^k_\alpha := \mathcal L^1_{\alpha^1} \circ \ldots \circ \mathcal L^k_{\alpha^k}, \qquad \wt{\mathcal R}^k_\alpha := \mathcal R^k_{\alpha^k}  \circ \ldots \circ \mathcal R^1_{\alpha^1}.
   $$
   Our aim is to estimate at each step $k$ the commutator $\bs \Gamma_k$ with the commutator $\bs \Gamma_{k+1}$, up to removing a small set where the new four left compositions $\mathcal L^{k+1}_{p}$, $p = 1, 2, 3, 4$ are not $e^{\lambda_i}$-Lipschitz. 

   Notice that for every $j\geq 1$, $\mathcal L^{j}_p$ is a composition of the flows of $X, Y$ for a total time of order $2^{-j}$, as it immediately follows from \eqref{eq:right}. Recall that $g = g_X + g_Y$ is the sum of the maximal functions of $DX, DY$, which are in $\mathbf L^p_{loc}$.
   By Lemma \ref{lemma:Lusinlip}, for all $j\geq 0$ and $p \in \{1, 2, 3, 4\}$, there exist sets $\mathcal B^j_p$ such that
$$
|\mathcal B^j_p|\,  \lesssim \, \frac{1}{\lambda_j 2^j} \int_{\{g \geq \lambda_j 2^j/(2T)\}} \, g(x) \, \dif x
$$
and
   \begin{equation}
      \mathcal L^j_{p}\quad \text{is $e^{c_d\lambda_j}$-Lipschitz on $\mathbb R^d \setminus \mathcal B^j_p$.}
   \end{equation}
   Define the bad set:
$$
\begin{aligned}
\mathcal B_j & := \bigcup_{\alpha \in \mathcal I^j} \bigcup_{p=1}^4 \Bigg\{x \in \mathbb R^d \; \Big| \; A_j \circ B_j \circ \wt{\mathcal R}^j_\alpha(x) \in \mathcal B^j_p \Bigg\} \\
&  \qquad  \qquad \bigcup_{\alpha \in \mathcal I^j} \bigcup_{p=1}^4 \Bigg\{x \in \mathbb R^d \; \Big| \; B_j \circ A_j \circ \wt{\mathcal R}^j_\alpha(x) \in \mathcal B^j_p \Bigg\}.
\end{aligned}
$$
Therefore we estimate
\begin{equation}
    |\mathcal B_j| \leq 2 \cdot 4 \cdot 4^j  \cdot \frac{1}{\lambda_j 2^j} \int_{\{g \geq \lambda_j 2^j/(2T)\}} \, g(x) \, \dif x.
\end{equation}
The \emph{good set} at stage $k$ is then 
$$
\mathcal G_k :=  \mathbb R^d \, \setminus \, \bigcup_{j=1}^k \mathcal B_j.
$$
Therefore, in the good region $\mathcal G_k$ we have that every appearing left factor $\mathcal L^j_p$ is $e^{\lambda_j}$ Lipschitz on the points given by the right composition.

\vspace{0.3cm}
 \noindent {\bf 2.} We want to show that for every ${k^\prime} \leq k$, there holds
\begin{equation}\label{eq:induction}
\begin{aligned}
    \int_{\mathcal G_k\cap B_r} |\bs \Gamma_0(x)| \dif x & \leq   \sum_{\alpha \in \mathcal I^{k^\prime} } \Bigg(\prod_{j=1}^{k^\prime}  e^{c_d\lambda_j}\Bigg)\int_{\mathcal G_k\cap B_{r}}  |\bs \Gamma_{k^\prime} \circ \wt{\mathcal R}_\alpha^{k^\prime} | \dif x.
\end{aligned}
\end{equation}
This can be done, for example, by an easy induction. For $k^{\prime} = 0$, the result is trivial (here for $k = 0$, we set $\mathcal R^0_\alpha =  \mr{Id}$).

Assume to have proved \eqref{eq:induction} up to $k^\prime \in \mathbb N$, where $1\leq  k^\prime < k$. Decomposing $\bs \Gamma_{k^\prime}$ as in \eqref{eq:jdec}, we obtain 
\begin{equation}
\begin{aligned}
   \sum_{\alpha \in \mathcal I^{k^\prime} } & \Bigg(\prod_{j=1}^{k^\prime}  e^{c_d\lambda_j}\Bigg)\int_{\mathcal G_k\cap B_r}  |\bs \Gamma_{k^\prime} \circ \wt{\mathcal R}_\alpha^{k^\prime} | \dif x\\
    & \leq \sum_{p=1}^4 \sum_{\alpha \in \mathcal I^{k^\prime} } \Bigg(\prod_{j=1}^{k^\prime}  e^{c_d\lambda_j}\Bigg)\int_{\mathcal G_k\cap B_r}  |\mathcal L^{k^\prime+1}_p \circ \bs \Gamma_{k^\prime+1} \circ \mathcal R^{k^\prime+1}_p \circ \wt{\mathcal R}_\alpha^{k^\prime} | \dif x\\
    & =  \sum_{\alpha \in \mathcal I^{k^\prime+1} } \Bigg(\prod_{j=1}^{k^\prime}  e^{c_d\lambda_j}\Bigg)\int_{\mathcal G_k\cap B_r}  |\mathcal L^{k^\prime+1}_{ \alpha^{k^\prime+1}} \circ \bs \Gamma_{k^\prime+1} \circ \wt{\mathcal R}^{k^\prime+1}_\alpha | \dif x \\
    & \leq \sum_{\alpha \in \mathcal I^{k^\prime+1} } \Bigg(\prod_{j=1}^{k^\prime+1}  e^{c_d\lambda_j}\Bigg)\int_{\mathcal G_k\cap B_r}  |\bs \Gamma_{k^\prime+1} \circ \wt{\mathcal R}^{k^\prime+1}_\alpha | \dif x.
    \end{aligned}
\end{equation}
Here in the equality of the third line we used the fact that if $\alpha \in \mathcal I^{k^\prime}$ and $p = \alpha^{k^\prime+1}$, then $\alpha^\prime := (\alpha, p) \in \mathcal I^{k^\prime+1}$, and by definition
$$
\mathcal R^{k^\prime+1}_p \circ\wt{\mathcal R}_\alpha^{k^\prime} = \wt{\mathcal R}_{\alpha^\prime}^{k^\prime+1} 
$$
In the last inequality, we used the fact that $\mathcal L^{k^\prime+1}_{\alpha^{k^\prime+1}}$ is calculated only on pair of points where it is $e^{\lambda_{k^\prime+1}}$-Lipschitz thanks to the definition of the good set $\mathcal G_k$.
This proves the inductive step, so that \eqref{eq:induction} holds in particular for $k^\prime = k$, for every $k \geq 1$.

Next, let $k \in \mathbb N$, apply \eqref{eq:induction} with $k^\prime = k$, and write 
\begin{equation}
    \begin{aligned}
        \int_{B_r} |\bs \Gamma_0(x)| \dif x \, & \lesssim \, \sum_{j=1}^k |\mathcal B_j| +\sum_{\alpha \in \mathcal I^{k} } \Bigg(\prod_{j=1}^{k}  e^{c_d\lambda_j}\Bigg)\int_{\mathcal G_k\cap B_r}  |\bs \Gamma_{k} \circ \wt{\mathcal R}_\alpha^{k} | \dif x \\
        & \, \lesssim\,  \sum_{j=1}^k \frac{2^j}{\lambda_j} \int_{\{g \geq \lambda_j 2^j /(2T)\}} \, g(x) \, \dif x \\
        & \qquad + 4^k \Bigg(\prod_{j=1}^{k}  e^{c_d\lambda_j}\Bigg)\int_{B_{R_k}}  |\bs \Gamma_{k}(x) | \dif x.
    \end{aligned}
\end{equation}
Here in the last step we used that every right factor $\wt{\mathcal R}_\alpha^k$ has uniformly bounded compression, since it is a composition of flows of time at most $2T$, and that the points $\wt{\mathcal R}_\alpha^k$ are contained in the ball $B_{R_k}$.
\end{proof}

We are now in the position of proving the main theorem for vector fields in $W^{1,2}$.
\begin{thm}
Let $X, Y \in W^{1,2}(\mathbb R^d) \cap \mathbf L^\infty$ with bounded divergence, and assume that $[X, Y] = 0$. Then for all $s, t$, there holds
\begin{equation}
    \Phi^t_X \circ \Phi^s_Y(x) - \Phi^s_Y \circ \Phi^t_X(x) = 0 \qquad \text{for a.e. $x \in \mathbb R^d$}.
\end{equation}
\end{thm}
\begin{proof}
   Let $ r >0$ be fixed and take $0 \leq s, t \leq T = 1$.  From Proposition \ref{prop:key}, we deduce
   \begin{equation}\label{eq:key-estimate}
\begin{aligned}
&\int_{B_r}
\Big|\Phi_X^t\circ\Phi_Y^s(x)
-\Phi_Y^s\circ\Phi_X^t(x)\Big|\dif x\\
&\quad\lesssim
\sum_{j=1}^{k}
\frac{2^j}{\lambda_j}
\int_{ \{g>\,2^j\lambda_j/(2T)\}}g(x)\dif x\\
&\qquad+
4^k\left(\prod_{j=1}^{k}e^{c_d\lambda_j}\right)
\int_{B_{R_k}}
\Big|
\Phi_X^{\frac{t}{2^k}}
\circ\Phi_Y^{\frac{s}{2^k}}(x)
-
\Phi_Y^{\frac{s}{2^k}}
\circ\Phi_X^{\frac{t}{2^k}}(x)
\Big|\dif x.
\end{aligned}
\end{equation}
Let $\eta > 0$ be fixed but sufficiently small, and $M > 0$.  Define
 $$
 \lambda_j  := \eta + M2^{-j}
 $$
 where $\eta$ is fixed and sufficiently small to be chosen later.
 The first term is then bounded by 
 \begin{equation}\label{eq:first}
     \begin{aligned}
        \sum_{j=1}^{k}
\frac{2^j}{ \lambda_j}
\int_{ \{g>\,2^j\lambda_j/(2T)\}} g(x)\dif x & = \int_{\{g \geq M/(2T)\}} g(x) \sum_{j=1}^{k}
\frac{2^j}{\lambda_j} \mathbf 1_{\{g>\,2^j\lambda_j/(2T)\}}(x) \dif x \\
& \leq \frac{1}{\eta}  \int_{\{g \geq M/(2T)\}} g(x) \sum_{j=1}^{\lfloor\log_2 (2Tg/\eta)\rfloor} 2^j \dif x\\
& \leq \frac{4T}{\eta^2} \int_{\{g \geq M/(2T)\}} g^2(x) \dif x.
     \end{aligned}
 \end{equation}
 For the second term, we use Proposition \ref{prop:2}, so that we obtain 
\begin{equation}
\begin{aligned}
      4^k\left(\prod_{j=1}^{k}e^{c_d\lambda_j}\right)
& \int_{B_{R_k}}
\Big|
\Phi_X^{\frac{t}{2^k}}
\circ\Phi_Y^{\frac{s}{2^k}}(x)
-
\Phi_Y^{\frac{s}{2^k}}
\circ\Phi_X^{\frac{t}{2^k}}(x)
\Big|\dif x \\
      & \lesssim 4^k  e^{c_d(k\eta +2M)} 2^{-k} |\{g \geq 2^k/T\} | \\
      & \lesssim  e^{2c_d M} e^{c_dk\eta}2^{-k}  \int_{\{g \geq 2^k/T\} } g^2(x) \dif x
\end{aligned}
\end{equation}
If $\eta$ is sufficiently small ($\eta \leq \log(2)/c_d$), the right hand side goes to zero as $k \to +\infty$ and $M$ is fixed. Then, we let $M$ go to infinity, so that the right hand side of \eqref{eq:first} goes to zero as well. This proves the result.
\end{proof}

For $p > 2$, there is a clean quantitative version using a different choice of the sequence $\lambda_i.$
\begin{thm}
Let $X, Y \in W^{1,p}(\mathbb R^d) \cap \mathbf L^\infty$ with $2 < p < \infty$, with bounded divergence. Then for all $0 \leq s, t \leq T$, there holds
\begin{equation}
 \int_{B_r} \left|\Phi^t_X \circ \Phi^s_Y(x) - \Phi^s_Y \circ \Phi^t_X(x)\right| \dif x \lesssim  \left[ \log \left( 1 + \frac{1}{\|[X, Y] \|_{\mathbf L^1(B_R)}}\right) \right]^{-p}.
    \end{equation}
    where $R$ is as in Proposition \ref{prop:key}.
\end{thm}

\begin{proof}
We apply the same argument as in the proof of the theorem above. This time, choose 
$$
\lambda_j = \lambda 2^{-\alpha j}, \qquad 
\alpha = \frac{p-2}{2p}.
$$
From \eqref{eq:key-estimate}, we can now estimate the first term in the right hand side as 
 \begin{equation}\label{eq:first1}
     \begin{aligned}
        \sum_{j=1}^{k}
\frac{2^j}{ \lambda_j}
\int_{\{g>\,2^j\lambda_j/(2T)\}} g(x)\dif x & \leq \frac{1}{\lambda^p} \left( \sum_{j=1}^k 2^{j\left(1-\frac{p}{2}\right)}\right) \, \int_{\mathbb R^d} g^p(x) \dif x \\
& \leq \frac{C_p}{\lambda^p}\, \int_{\mathbb R^d} g^p(x) \dif x
     \end{aligned}
 \end{equation}
 where $C_p$ is a constant depending on $p$.
 The second term in \eqref{eq:key-estimate} is now estimated by Proposition \ref{prop:2} as 
\begin{equation}\label{eq:above}
\begin{aligned}
      4^k\left(\prod_{j=1}^{k}e^{c_d\lambda_j}\right)
& \int_{B_{R_k}}
\Big|
\Phi_X^{\frac{t}{2^k}}
\circ\Phi_Y^{\frac{s}{2^k}}(x)
-
\Phi_Y^{\frac{s}{2^k}}
\circ\Phi_X^{\frac{t}{2^k}}(x)
\Big|\dif x \\
   &  \lesssim 4^k e^{\wt C_p \lambda} \Big( 2^{-k} \, |\{g \geq 2^{k}\}|\, +  2^{-2k} \,\big\| [X,Y]\big\|_{\mathbf L^1(B_R)}\Big) \\
   & \lesssim e^{\wt C_p \lambda}  2^{-(p-1)k} \int_{\mathbb R^d} g^p(x) \dif x + e^{\wt C_p \lambda} \,\big\| [X,Y]\big\|_{\mathbf L^1(B_R)}.
\end{aligned}
\end{equation}
Then, send $k$ to infinity. Since $p > 2$ the first term in \eqref{eq:above} vanishes, and we get
\begin{equation}
    \int_{B_r}
\Big|\Phi_X^t\circ\Phi_Y^s(x)
-\Phi_Y^s\circ\Phi_X^t(x)\Big|\dif x \leq  \frac{C_p}{\lambda^p}\, \int_{\mathbb R^d} g^p(x) \dif x + e^{\wt C_p \lambda} \,\big\| [X,Y]\big\|_{\mathbf L^1(B_R)}.
\end{equation}
Finally, optimizing in $\lambda$ we obtain the result.
   \end{proof}

\section{Non-commutativity for $1 \leq p < 2$}
\begin{prop}
    Let $p \in [1, 2)$. Then there exist divergence free vector fields $X, Y \in W^{1,p}_{loc}(\mathbb R^2)$ with $[X,Y] = 0$ such that their regular Lagrangian flows do not commute.
\end{prop}

\begin{proof}
    The vector fields $X ,Y$ are constructed using a double covering of the plane. Consider the map $Q : \mathbb R^2 \to \mathbb R^2$ which doubles the angle. In polar coordinates, 
    $$
Q(\rho, \theta) := (\rho, 2\theta).
$$
The map $Q$ is clearly measure preserving for the Lebesgue measure, since it is in particular measure preserving on every circle $\{|x| = \rho\}$. Consider the vector fields $e_1, e_2$ in $\mathbb R^2$, constant and parallel to the coordinates axis, with unit norm. Clearly, $e_1, e_2$ are divergence free and $[e_1, e_2] = 0$. We define the pullbacks of the vector fields $e_1, e_2$ by
$$
X(x) := (D Q(x))^{-1}\,e_1, \qquad Y(x) := (D Q(x))^{-1}\,e_2.
$$
$Q$ is measure preserving and has Jacobian $\det(DQ(x)) = 2$ for all $x \neq 0$, therefore $X, Y$ are divergence free. Moreover, the pullback of vector fields naturally preserves the Lie brackets, therefore 
$$
[X, Y] = 0.
$$
In polar coordinates, one has
\begin{equation}\label{eq:branchfields}
 X=\cos(2\theta)e_r-\tfrac12\sin(2\theta)e_\theta,
 \qquad
 Y=\sin(2\theta)e_r+\tfrac12\cos(2\theta)e_\theta.
\end{equation}
It is therefore not difficult to show that 
\[
 |D\,X(x)|+|D\,Y(x)|\, \sim \,  \frac{1}{|x|}
\]
therefore $X, Y \in W^{1,p}_{loc}(\mathbb R^2)$ if and only if $p \in [1,2)$.
Consider a small neighborhood $U$ of the point $(-1,-1) \in \mathbb R^2_y$ and consider trajectories starting from $y \in U$, obtained by first following $e_1$ for time $2$ and $e_2$, also for time $2$, and conversely. Then the corresponding pulled back trajectories for $X, Y$ in the $x$-plane do not commute, because  they belong to different branches of $Q$ (see Figure \ref{fig:nocomm}). To see this, fix $y_0=(-1,-1)$ and choose any
$x_0\in Q^{-1}(y_0)$. Consider the two
polygonal paths in $\mathbb R^2_y$:
\[
\begin{aligned}
\gamma_1 &: (-1,-1)\to(1,-1)\to(1,1),\\
\gamma_2 &: (-1,-1)\to(-1,1)\to(1,1)
\end{aligned}
\]
each one made of two segments. Since both paths avoid the origin, they can be uniquely lifted through $Q$, starting from $x_0$. By construction, the endpoints of these lifts are respectively $\Phi_Y^2\circ\Phi_X^2(x_0)$ and $\Phi_X^2\circ\Phi_Y^2(x_0)$.
The closed curve obtained by concatenating $\gamma_1$ with the reverse of $\gamma_2$ winds once around the origin. Since $Q$ is a two-branch covering of $\mathbb R^2$, the two lifts have distinct
endpoints. Hence
\[
\Phi_Y^2\circ\Phi_X^2(x_0)
\neq
\Phi_X^2\circ\Phi_Y^2(x_0).
\]
The same argument applies to all starting points $y_0$ in a sufficiently small neighborhood of $(-1,-1)$, since the corresponding rectangles still enclose the origin. Thus the two compositions differ on an open set of positive Lebesgue measure.
\end{proof}

\end{document}